\documentclass{amsart}
\newcommand{\CO}[2]{[#1,#2)}

\usepackage{amssymb,amsmath,amsthm}
\usepackage[usenames]{color}	
\usepackage[mathscr]{eucal}  
\usepackage{amsfonts,amssymb,latexsym,amsmath,amsthm,epsfig,amscd}

\newtheorem{lemma}{Lemma}[section]
\newtheorem{cor}[lemma]{Corollary}

\newtheorem{theorem}[lemma]{Theorem}

\theoremstyle{definition}

\newtheorem{example}[lemma]{Example}
\newtheorem {observation}[lemma]{Observation}

\newcommand{\labt}[1]{\label{theorem:#1}}
\newcommand{\reft}[1]{Theorem~\ref{theorem:#1}}

\newcommand{\labeq}[1]{\label{eq:#1}}

\newcommand{\R}{\mathbb{R}}
\newcommand{\Z}{\mathbb{Z}}
\newcommand{\QQ}{\mathcal Q}

\newcommand{\N}{\mathbb{N} }

\newcommand{\D}{{\mathcal{D}}_{\infty}^Q(\alpha)}
\newcommand{\Qp}{q_1\cdots q_n}

\newcommand{\BB}{\ub}

\newcommand\HD[1]{\text{{\rm HD}}\left( #1 \right)}
\newcommand {\ignore}[1]{}

\newcommand{\lp}{P_{Q,\alpha}}
\newcommand{\ub}{b_Q(\alpha)}

\newcommand{\floor}[1]{\left\lfloor #1 \right\rfloor}

\thanks{The first-named author was supported in
part by the Simons Foundation grant 245708. The second-named author was supported in part by the NSF grant DMS-0943870. The last-named author was supported in part by NSF Grant DMS 1361677.}

\author[L. Fishman]{Lior Fishman}
\author[B. Mance]{Bill Mance}
\author[D. Simmons]{David Simmons}
\author[M. Urba\'{n}ski]{Mariusz Urba\'{n}ski}

\IfFileExists{marvosym.sty}{
\usepackage{marvosym} \newcommand{\at}{\MVAt}
}{
\newcommand{\at}{@}
}

\address{University of North Texas, Department of Mathematics, 1155 Union Circle \#311430, Denton, TX 76203-5017, USA}
\email[L. Fishman]{lior.fishman\at unt.edu}
\email[B. Mance]{mance\at unt.edu}
\email[M. Urba\'{n}ski]{urbanski\at unt.edu}

\address{Ohio State University, Department of Mathematics, 231 W. 18th Avenue Columbus, OH 43210-1174, USA}
\email[D. Simmons]{simmons.465\at osu.edu}

\begin{document}

\title{Shrinking targets for non-autonomous dynamical systems corresponding to Cantor series expansions}

\begin{abstract}
We provide a closed formula of Bowen type for the Hausdorff dimension of a very general shrinking target scheme generated by the non-autonomous dynamical system on the interval $\CO 01$, viewed as $\R/\Z$, corresponding to a given method of Cantor series expansion. We also examine a wide class of examples utilizing our theorem. In particular, we provide a Diophantine approximation interpretation of our scheme.
\end{abstract}

\maketitle

\section{Introduction}

Recall that in the framework of autonomous dynamical systems, the evolution process of the system is determined by continuously iterating a fixed map. In contrast, the system is said to be non-autonomous if at each stage of the iteration process we allow the action of a possibly different map.
In this paper we will be concerned with a shrinking target problem in the context of the
non-autonomous dynamical system generated by any sequence $Q=(q_n)\in\N _{\geq 2}^\N$. Given such a sequence, the maps $T_{Q,n}:\R/\Z\to\R/\Z$ and $T_Q^n:\R/\Z\to\R/\Z$ are defined for any $n\in\N$ as follows:
\begin{equation}\label{the1st}
\begin{split}
T_{Q,n}(x) &= q_n x \pmod{1}\\
T_Q^n(x) &= T_{Q,n} \circ \cdots \circ T_{Q,1}(x) = q_n\cdot q_{n-1}\cdots q_1 x \pmod{1}.
\end{split}
\end{equation}
This system was first introduced and investigated in an implicit manner by G. Cantor in \cite{Cantor}, where he considered what is now called the {\sl{$Q$-Cantor series expansion}} of a real number $x$, i.e. the (unique) expansion of the form
\[
x=\omega_0+\sum_{j=1}^{\infty} \frac {\omega_j} {q_1 q_2 \cdots q_j}
\]
where $\omega_0=\floor{x}$, $\omega_j$ is in $\{0,1,\cdots,q_j-1\}$ for $j\geq 1$, and $\omega_j \neq q_j-1$ for infinitely many $j$. The relation between this definition and the non-autonomous dynamical system \eqref{the1st} is that for every $n\geq 1$, $\omega_n = \lfloor q_n T_Q^{n - 1}(x)\rfloor$, where $T_Q^{n - 1}(x)$ denotes the representative of $T_Q^{n - 1}(x)$ in $\CO 01$.

The study of normal numbers and various statistical properties of real numbers with respect to large classes of Cantor series expansions dates back to P. Erd\H{o}s and A. R\'{e}nyi \cite{ErdosRenyiConvergent, ErdosRenyiFurther}, and was continued by A. R\'{e}nyi \cite{RenyiProbability, Renyi, RenyiSurvey} and P. Tur\'{a}n \cite{Turan}.
Later, certain aspects of this area were extensively studied by many authors, notably by J. Galambos \cite{Galambos2}, T. \u{S}al\'{a}t (e.g. \cite{Salat2}), and F. Schweiger \cite{SchweigerCantor}. Most recently, the second-named author and  his collaborators continued to develop this area of research in \cite{AireyMance1, AireyManceVandehey, ppq1},
where the primary results concern various concepts of normality, relations between them, and the Hausdorff dimensions of appropriate significant sets. 

In this paper we follow a different approach. We are prompted by the (relatively) recent activity focused on determining the Hausdorff dimension of the set of points ``hitting'' some shrinking target infinitely often. As with research regarding Cantor expansions, one may trace this approach to the pioneering work of  A. Besicovi\u c \cite{Be} and V. Jarn\'{\i}k \cite{Ja} with respect to continued fraction expansions. Others investigate conformal dynamics, Kleinian groups, and conformal iterated function systems. We list here a (by no means exhaustive) selection for the reader's convenience: \cite{BL, BW, HV1, HV2, Reeve, St, SU, Urbanski}. We would however like to emphasize that unlike all the papers mentioned above, we work in the context of a non-autonomous dynamical system, namely the one defined in \eqref{the1st}.

The shrinking target scheme considered in the present paper is quite general, at least in the context of sequences. Let $\alpha =(\alpha_i)_{i=1}^\infty$ be a sequence of nonnegative real numbers, and for each $n\geq 1$ let 
\[
\alpha (n):=\sum_{i=1}^{n}\alpha_i.
\]  
Let
\[
\D = \{ x \in X := \R/\Z: \|T_Q^n(x)\|\leq e^ {-\alpha(n)} \text{ for infinitely many } n\},
\]
where $\|\cdot\|$ denotes distance to the nearest integer. We would like to bring the reader's attention to the fact that the set $\D$ has a precise Diophantine approximation interpretation. A general scheme of Diophantine analysis has been laid down by three of the authors in \cite{FSU}. In the setting considered here, for every integer $n\ge 1$, denote by $\QQ_n$ the set of all $Q$-adic rationals of order $n$ in $X$, i.e. the set of all numbers of the form
\[
\sum_{j=1}^n\frac{\omega_j}{q_1\cdots q_j},
\]
where $\omega_j\in \{0,1,\ldots,q_j-1\}$. For every $Q$-adic rational number $w\in\QQ := \bigcup_{n\geq 1} \QQ_n$, let $H_Q(w)$ denote the least integer $n\ge 1$ such that $w\in \QQ_n$. We can interpret $H_Q(w)$ as the {\sl{height}} (induced by the sequence $Q$) of the number $w$. The triple $(X,\QQ, H_Q)$ then forms a {\sl{Diophantine space}}, i.e., a complete metric space, a dense set, and a function measuring the ``complexity'' of elements of this dense set.  Define the (approximation) function $\psi:\N \to (0,\infty)$ as follows:
\[
\psi(n)= \frac{\exp(-\alpha(n))}{q_1\cdots q_n}\cdot
\]
In Diophantine approximation terminology, a point $x\in X$ is called \emph{$\psi$-approximable} (relative to the Diophantine space $(X,\QQ,H_Q)$) if there exists a sequence $(w_k)_1^\infty$ in $\QQ$ converging to $x$ with the property that $|x-w_k|\le \psi(H_Q(w_k)) \; \forall k\geq 1$. We observe then that $\D$ is precisely the set of $\psi$-approximable numbers, a set which is a standard object of study in Diophantine analysis.



Using the definitions above, and prompted by thermodynamic formalism, given $s\geq 0$, we define the {\sl{(upper) pressure function}} 
\begin{equation}
\label{upperpressure}
\lp(s)
:=\limsup_{n\to\infty}\frac{1}{n}[(1-s)\log(\Qp)-s\alpha(n)].
\end{equation}
Note that since $q_i\geq 2$ and $\alpha_i\geq 0$, the function $s\mapsto \lp(s)$ is strictly decreasing in its domain of finiteness.  We further define
\[
\ub := \sup\{t\geq 0 : P_{Q,\alpha}(t) > 0 \} = \inf \{t\geq 0 :P_{Q,\alpha}(t)<0\}.
\]


\begin{theorem}\labt{main}
For any sequence $\alpha =(\alpha_i)_{i=1}^\infty$ of nonnegative real numbers, 
\[
 \HD{\D}=  \ub = \limsup_{n\to \infty} \frac{\log(\Qp)}{\log(\Qp) + \alpha(n)}\cdot
\]
\end{theorem}

Here and in what follows, $\rm{HD}$ denotes Hausdorff dimension.

The general form of this theorem does not differ much from the ones obtained in the papers mentioned above. What is perhaps surprising is that it holds in such generality in the realm of a non-autonomous system. We note that it also covers the autonomous case of $q$-ary sequences, where $q\geq 2$ is an integer; simply consider the constant sequence $Q=(q)_1^\infty$ every term of which is equal to $q$. It also captures the cases of periodic and eventually periodic sequences $Q$, i.e. the ones that can be also approached with the methods of autonomous dynamical systems.

We prove \reft{main} in the next section and in Section~\ref{examples} we describe a number of classes of examples which illustrate its content and scope. 
 
\section{Proof of \reft{main}}


Fix $t > \BB$ arbitrary, so that $\lp(t) < 0$. For each $n\geq 1$, let $Q_n = q_1\cdots q_n$, and for each $j = 0,\ldots,Q_n - 1$, let
\[
\Delta_{n,j}(\alpha) := \left\{x\in X : \left\|x - \frac{j}{Q_n}\right\|\leq \frac{e^{-\alpha(n)}}{Q_n}\right\},
\]
so that
\[
\D = \limsup_{n\to\infty} \bigcup_{j = 0}^{Q_n - 1} \Delta_{n,j}(\alpha) = \bigcap_{N = 1}^\infty \bigcup_{n = N}^\infty \bigcup_{j = 0}^{Q_n - 1} \Delta_{n,j}(\alpha).
\]
We have
\[
|\Delta_{n,j}(\alpha)|= 2(Q_n )^{-1}e^{-\alpha(n)}
\] 
and thus for all sufficiently large $n\geq 1$,
\[
\sum_{j = 0}^{Q_n - 1} |\Delta_{n,j}(\alpha)|^t = Q_n (2(Q_n )^{-1}e^{-\alpha(n)})^t = 2^t Q_n^{1 - t} e^{-t\alpha(n)} \leq 2^t \exp\left(\frac12\lp(t) \cdot n\right).
\]
Thus for all sufficiently large $N\geq 1$,
\[
\sum_{n = N}^\infty \sum_{j = 0}^{Q_n - 1} |\Delta_{n,j}(\alpha)|^t \leq 2^t \sum_{n = N}^\infty \exp\left(\frac12\lp(t) \cdot n\right)
\]
and so
\[
H^t(\D) \leq 2^t \lim_{N\to\infty} \sum_{n = N}^\infty \exp\left(\frac12\lp(t) \cdot n\right) = 0,
\]
where $H^t$ denotes Hausdorff $t$-dimensional measure. Hence $\HD{\D} \leq t$, and since $t > \ub$ was arbitrary $\HD{\D} \leq \ub$.

We now prove that $\HD{\D}\geq \ub$. Fix $0\leq s < \ub$ and fix a sufficiently fast growing sequence $(n_i)_{i=1}^\infty$, to be
determined later in the proof, along which the lim sup in \eqref{upperpressure} is achieved, i.e. for which
\begin{equation}
\label{requirement}
\lim_{l\to\infty} \frac{1}{n_l}[(1-s)\log(Q_{n_l})-s\cdot \alpha(n_l)]=\lp(s)>0.
\end{equation}
For each $l\geq 1$ let
\[
S_l = \{\Delta_{n_l,j}(\alpha) : j = 0,\ldots, Q_{n_l} - 1\}.
\]
We assume that $n_1$ is chosen large enough so that $\alpha(n_1) > \log(2)$ (ignoring the case $\alpha(n)\leq \log(2)\;\forall n$ as trivial), so that for each $l\geq 1$, $S_l$ is a disjoint collection. Now we construct inductively a sequence of sets $(R_l)_{l=1}^\infty$ as follows. We start by letting $R_1 = S_1$. For the inductive step, suppose the set $R_l\subseteq S_l$ has been defined. Then we let $R_{l+1}$ be the set consisting of all elements $\Delta\in S_{l + 1}$ contained in some interval from the family $R_l$. We define the following nonempty compact set:
\[
K:= \bigcap_{l = 1}^\infty \bigcup(S_l) = \bigcap_{l=1}^\infty \bigcup_{\Delta\in S_l} \Delta.
\]
For every $\Delta\in S_l$ let
\[
R_{l+1}(\Delta):=\{\Gamma\in R_{l+1}: \Gamma\subseteq\Delta\}.
\]
Then
\[
R_{l+1}=\bigcup_{\Delta\in S_l}R_{l+1}(\Delta)
\]
and for all $\Delta\in S_l$,
\begin{equation}\label{cardinalityMU}
\begin{aligned}
\#(R_{l+1}(\Delta))
&\ge \frac{|\Delta|}{\left(Q_{n_{l+1}}\right)^{-1}}-2 
=\frac{(Q_{n_l})^{-1}e^{-\alpha(n_l)}}{\left(Q_{n_l} q_{{n_l}+1}\cdots q_{n_{l+1}}\right)^{-1}}-2 \\
&=q_{n_l+1}\cdots q_{n_{l+1}}e^{-\alpha(n_l)}-2
\ge\frac12 q_{n_l+1}\cdots q_{n_{l+1}}e^{-\alpha(n_l)},
\end{aligned}
\end{equation}
where the last inequality holds provided that $n_{l+1}\ge 1$  is large enough.
We now shall recursively define maps $m_l:R_l\to[0,1]$, $l\ge 1$, as follows. Let $m_1(\Delta) = 1/\#(R_1)$ for all $\Delta\in R_1$. Proceeding inductively, fix $\Delta\in R_l$ and for every $\Gamma\in R_{l+1}$, set
\[
m_{l+1}(\Gamma)
:=\frac{m_{l}(\Delta)}{\#(R_{l+1}(\Delta))}
\]
Then it is easy to see (e.g. by choosing arbitrary measures extending the functions $m_l$ and then taking a weak limit) that there exists a Borel probability measure $m$ supported on $K$ such that $m(\Delta) = m_l(\Delta)$ for all $l\geq 1$ and $\Delta\in S_l$. Now \eqref{cardinalityMU} shows that
\[
m_{l+1}(\Gamma) \le 2m_{l}(\Delta)e^{\alpha(n_l)}\cdot(q_{n_l+1}\cdots q_{n_{l+1}})^{-1}
\]
and iterating this estimate gives
\begin{equation}\label{cylindermeasureestimate}
m(\Delta)
\le 2^{l-1}\exp\left(\alpha(n_{1})+\alpha(n_{2})+\ldots+\alpha(n_{l-1})\right)
    \cdot(Q_{n_l})^{-1}.
\end{equation}
Let $x\in K$ be arbitrary. We want to show that there exists some $C>0$ independent of $x$  such that for every $r>0$
\begin{equation}\label{108052014}
m(B(x,r))\leq C\cdot r^s.
\end{equation}
Since $m$ is a probability measure, it is of course enough to show this for all $r>0$ small enough. Fix $r\in (0, e^{-\alpha(n_1)}(Q_{n_1})^{-1})$ and then  let  $l$ be the largest integer such that 
\begin{equation}\label{condition}
e^{-\alpha(n_l)}\cdot(Q_{n_l})^{-1}\geq r.
\end{equation}
By our choice of $r$, we have that $l\ge 1$. Cover $B(x,r)$ by a union of intervals of the form $\left[\frac{j - 1/2}{Q_n}, \frac{j + 1/2}{Q_n}\right]$, $j = 0,\ldots,Q_n - 1$. We can do it by taking no more than
\[
\frac{r}{(Q_{n_{l + 1}})^{-1}}+2\le 2rQ_{n_{l + 1}}
\] 
such intervals. But then we can also cover $K\cap B(x,r)$ by at most $2rQ_{n_{l + 1}}$ intervals of $R_{l + 1}$. Invoking \eqref{cylindermeasureestimate} and the fact that the measure $m$ is supported on $K$, we therefore get that
\[
\begin{aligned}
m(B(x,r))
&\le 2rQ_{n_{l + 1}}2^{l}\exp\left(\alpha(n_{1})+\alpha(n_{2})+\ldots+\alpha(n_{l})\right)
    \cdot Q_{n_{l + 1}}^{-1} \\
&=2^{l+1}r\exp\left(\alpha(n_{1})+\alpha(n_{2})+\ldots+\alpha(n_{l})\right).
\end{aligned}
\]
Hence, in order to show that \eqref{108052014} holds, it is enough to check that
\[
2^{l+1}\exp\left(\alpha(n_{1})+\alpha(n_{2})+\ldots+\alpha(n_{l})\right)
\leq C\cdot r^{s-1}.
\] 
Because of (\ref{condition}), and since $s < \ub \leq 1$, it is enough to show that
\[
2^{l+1}\exp\left(\alpha(n_{1})+\alpha(n_{2})+\ldots+\alpha(n_{l})\right)
\leq C(Q_{n_l})^{1-s}\exp\left((1-s)\alpha(n_{l})\right).
\]
Equivalently, 
\[
(l+1)\log 2+\alpha(n_{2})+\ldots+\alpha(n_{l-1}) \leq \log C + (1-s)\log(Q_{n_l})-s\cdot \alpha(n_{l}).
\]
But, because of our choice of $s$, we have that $\lp(s)>0$,  so by \eqref{requirement} we have
\[
(1-s)\log(Q_{n_l})-s\cdot \alpha(n_l)
\geq \frac{1}{2}\lp(s)\cdot n_l
\] 
for all large enough $l$. Thus, it suffices to verify that 
\[
\frac{1}{2}\lp(s)\cdot n_l
\ge (l+1)\log 2+\alpha(n_{1})+\alpha(n_{2})+\ldots+\alpha(n_{l-1}) - \log C.
\]
This can be done defining $n_l$ inductively to be large enough depending on $n_1, \dots ,n_{l-1} $ and on $l$. Note that this choice does not conflict with the requirement \eqref{requirement}. This completes the proof that $\HD{\D} = \ub$.

To finish the proof we need to show that $\ub = \delta := \limsup_{n\to\infty} Q_n/(Q_n + \alpha(n))$. If $s > t > \delta$, then for all $n$ sufficiently large we have
\begin{align*}
\frac{Q_n}{Q_n + \alpha(n)} &\leq t, &
\alpha(n) &\geq \frac{1 - t}{t}\log(Q_n)
\end{align*}
and thus
\begin{align*}
\lp(s) &\leq \limsup_{n\to\infty} \frac{1}{n}\left[(1-s)\log(Q_n)-s\frac{1 - t}{t}\log(Q_n)\right]\\
&= s\left(\frac{1 - s}{s} - \frac{1 - t}{t}\right)\liminf_{n\to\infty} \frac{\log(Q_n)}{n}
\leq s\left(\frac{1 - s}{s} - \frac{1 - t}{t}\right) \log 2 < 0
\end{align*}
and so $s > \ub$. Conversely, if $s < t < \delta$, then for infinitely many $n$ we have
\begin{align*}
\frac{Q_n}{Q_n + \alpha(n)} &\geq t, &
\alpha(n) &\leq \frac{1 - t}{t}\log(Q_n)
\end{align*}
and thus
\begin{align*}
\lp(s) &\geq \limsup_{n\to\infty} \frac{1}{n}\left[(1-s)\log(Q_n)-s\frac{1 - t}{t}\log(Q_n)\right]\\
&= s\left(\frac{1 - s}{s} - \frac{1 - t}{t}\right)\limsup_{n\to\infty} \frac{\log(Q_n)}{n}
\geq s\left(\frac{1 - s}{s} - \frac{1 - t}{t}\right) \log 2 > 0
\end{align*}
and so $s < \ub$.

\section{Examples}\label{examples}

In this section we consider a few classes of examples. Having proved \reft{main} our task reduces to examining sequences $(q_n)_{n=1}^\infty$ with appropriate arithmetical properties. Our examples show how to get a very explicit closed formula for the Hausdorff dimension of ${\D}$ for many classes of sequences $(q_n)_{n=1}^\infty$. 

\

\ignore{

We conclude with an example for which the lower and the upper Bowen's parameters $\lb$ and $\ub$ differ and the Hausdorff dimension $\HD{\D}$ cannot be determined by our theorem. 


}

Our first example follows directly from \reft{main}.

\begin{example}
Let $G(a_1,\cdots,a_n)$ denote the geometric mean of the positive real numbers $a_1,\cdots,a_n$.
Suppose that $Q$ is eventually periodic.  That is, we can write $Q$ in the form
\[
(d_1, d_2, \cdots, d_k, \overline{p_1, p_2, \cdots, p_m}).
\]
Let $\alpha_n = c > 0$ for all $n$. 
\end{example}

\noindent A short calculation shows that
\[
\HD{\D}=\frac {\log G(p_1,\cdots,p_m)} {\log G(p_1,\cdots,p_m)+c}.
\]
In particular, if $Q=(2,3,2,3,\cdots)$, then $\HD{\D}=\frac {\log \sqrt{6}} {\log \sqrt{6}+c}$.

 \noindent As mentioned it in the introduction, one could obtain this results by the methods of autonomous dynamics. 
This is particularly transparent in the case when $q_n=b$ for all $n$.

We observe that if
\[
L = \liminf_{n\to\infty} \frac{\alpha(n)}{\log(Q_n)},
\]
then \reft{main} says that
\begin{equation}
\label{main2}
\HD{\D} = \frac{1}{1 + L},
\end{equation}
with the convention that the right hand side is $0$ if $L = \infty$. On the other hand, we have the following simple observation:

\begin{observation}

\label{obs}
Let $L\in\R\cup\{\pm\infty\}$ and let $(a_n)_{n=1}^\infty$ and $(b_n)_{n=1}^\infty$ be two sequences of positive real numbers such that
\[
\sum_{n=1}^{\infty} b_n=\infty \hbox{ and } \lim_{n \to \infty} \frac {a_n} {b_n}=L.
\]
Then
\[
\lim_{n \to \infty} \frac {a_1+a_2+\ldots+a_n} {b_1+b_2+\ldots+b_n}=L.
\]
\end{observation}

Combining Observation \ref{obs} with \eqref{main2} yields the following:

\begin{cor}\label{cor}
Suppose that the limit
\begin{equation}\labeq{thm2L}
L := \lim_{n \to \infty} \frac {\alpha_n} {\log(q_n)}
\end{equation}
exists in $[0,\infty]$. Then
\[
\HD{\D}=\frac {1} {1+L},
\]
with the convention that the right hand side is $0$ if $L = \infty$.
\end{cor}

The next two examples follow directly from Corollary \ref{cor}.

\begin{example}
Suppose that $Q$ is a sequence such that $\lim_{n \to \infty} q_n / n^k \in (0, \infty)$ for $k>0$ and $\alpha_n=c \log n$ for $c>0$.  Then
\[
\HD{\D}=\frac {k} {k+c}.
\]
For example, if $q_n=n+1$, then $\HD{\D}=\frac {1} {1+c}$. As another example, if 
\[
q_n=2+\floor{\frac {\sqrt{n+\sqrt{n}\cos n}} {\sqrt[3]{n}}},
\]
then $\HD{\D}=\frac {1/6} {1/6+c}$.
\end{example}

\begin{example}
Suppose that $Q$ is a sequence such that $\lim_{n \to \infty} q_n / b^n \in (0,\infty)$ for $b>1$ and $\alpha_n=c n$ for $c>0$.  Then
\[
\HD{\D}=\frac {\log b} {\log b+c}.
\]
For example, if $q_n=2^n$, then $\HD{\D}=\frac {\log 2} {\log 2+c}$.
\end{example}

%
%
%

\end{document}